# Chow ring structure made simple

Sergey I. Nikolenko*     Nikita S. Semenov[†]

## Abstract

We show how to translate the task of computing the multiplicative structure of a Chow ring of a projective homogeneous variety into an easily understandable combinatorial task of calculating in the corresponding polynomial ring. The algorithms are also presented as a Maple package `ChowMaple` [Chow06].

Then we proceed to compute the multiplicative structure of the Chow rings for projective homogeneous varieties $E_6/P_1$, $E_7/P_7$, and $E_8/P_8$. For $E_6$, our results coincide with the results of Iliev and Manivel [IM05], who employed nontrivial geometric considerations to get the same results; our approach greatly simplifies the computations.

## 1   Introduction

The problem of calculating Chow rings has been among the most interesting problems of algebraic geometry for a long time. The Chow ring $CH(X)$ has proven to be a very powerful invariant, distinguishing nearly all known varieties. In [NSZ05], we described a way to calculate products in the Chow ring of $F_4$-varieties by using the Giambelli formula; for the applications listed there, namely calculating $CH(F_4/P_1)$ and $CH(F_4/P_4)$, the Giambelli formula sufficed. However, we had to make the algorithms more efficient in order to carry out the calculations for other projective homogeneous varieties.

In this paper, we describe how we have calculated the multiplicative structure of $CH(E_6/P_1)$, $CH(E_7/P_7)$, and $CH(E_8/P_8)$, which are probably the most interesting exceptional cases. Section 2 shows the connections between Chow rings and Hasse diagrams of the corresponding algebraic groups and describes the basics of our method, which were already developed in [NSZ05]. In Section 3 we discuss the improved algorithms in detail, while in Section 4 we briefly show the basics of our Maple package. Section 5 contains the main results of our work; it shows the multiplicative structure of the Chow rings.

One should mention the paper of Duan [Du05] devoted to another algorithm performing the same task. However, we believe that results concerning $E_8/P_8$ have never been published before.

*St.-Petersburg Department of the Steklov Institute of Mathematics. Supported by the Science School grant NSh-8464.2006.1.

[†]Universität Bielefeld. Supported partially by DAAD, INTAS.

# 2 Hasse diagrams and Chow rings

Let us fix a field $k$. We freely denote by $\Phi/P$ the variety $G/P$, where $G$ is a split simple algebraic group of type $\Phi$ over $k$, $P$ being its parabolic subgroup.

To each projective homogeneous variety $X$ we may associate an oriented labeled graph $\mathcal{H}$ called Hasse diagram. It is known that the ring structure of $\mathrm{CH}(X)$ is determined by $\mathcal{H}$. In the present section we remind several facts concerning relations between Hasse diagrams and Chow rings. For detailed explanations of these relations see [De74], [Hi82a] and [Ko91].

**2.1.** Let $G$ be a split simple algebraic group defined over a field $k$. We fix a maximal split torus $T$ in $G$ and a Borel subgroup $B$ of $G$ containing $T$ and defined over $k$. We denote by $\Phi$ the root system of $G$, by $\Pi = \{\alpha_1, \ldots, \alpha_{\mathrm{rk}G}\}$ the set of simple roots of $\Phi$ corresponding to $B$, by $W$ the Weyl group, and by $S$ the corresponding set of fundamental reflections.

Let $P = P_\Theta$ be a (standard) parabolic subgroup corresponding to a subset $\Theta \subset \Pi$, i.e., $P = BW_\Theta B$, where $W_\Theta = \langle s_\theta, \theta \in \Theta \rangle$. Denote

$$W^\Theta = \{w \in W \mid \forall\, s \in \Theta \quad l(ws) = l(w) + 1\},$$

where $l$ is the length function. The pairing

$$W^\Theta \times W_\Theta \to W \qquad (w, v) \mapsto wv$$

is a bijection and $l(wv) = l(w) + l(v)$. It is easy to see that $W^\Theta$ consists of all representatives in the cosets $W/W_\Theta$ which have minimal length. Sometimes it is also convenient to consider the set of all representatives of maximal length. We shall denote this set as $^\Theta W$. Observe that there is a bijection $W^\Theta \to {}^\Theta W$ given by $v \mapsto vw_\theta$, where $w_\theta$ is the longest element of $W_\Theta$. The longest element of $W^\Theta$ corresponds to the longest element $w_0$ of the Weyl group.

As $P_i$ we denote $P_{\Pi \setminus \{\alpha_i\}}$.

**2.2.** To a subset $\Theta$ of the finite set $\Pi$ we associate an oriented labeled graph, which we call a Hasse diagram and denote by $\mathcal{H}_W(\Theta)$. This graph is constructed as follows. The vertices of this graph are the elements of $W^\Theta$. There is an edge from a vertex $w$ to a vertex $w'$ labelled with $i$ if and only if $l(w) < l(w')$ and $w' = s_i w$. Observe that the diagram $\mathcal{H}_W(\emptyset)$ coincides with the Cayley graph associated to the pair $(W, S)$.

**2.3.** Now consider the Chow ring of a projective homogeneous variety $G/P_\Theta$. It is well known that $\mathrm{CH}(G/P_\Theta)$ is a free abelian group with a basis given by varieties $[X_w]$ that correspond to the vertices $w$ of the Hasse diagram $\mathcal{H}_W(\Theta)$. The degree of the basis element $[X_w]$ corresponds to the minimal number of edges needed to connect the respective vertex $w$ with $w_\theta$ (which is the longest word). The multiplicative structure of $\mathrm{CH}(G/P_\Theta)$ depends only on the root system of $G$ and the diagram $\mathcal{H}_W(\Theta)$.

By definition one immediately obtains

**2.4 Lemma.** *The contravariant functor* $\mathrm{CH}\colon \Theta \mapsto \mathrm{CH}(G/P_\Theta)$ *factors through the category of Hasse diagrams* $\mathcal{H}_W$, *i.e., the pull-back (ring inclusion)*

$$\mathrm{CH}(G/P_{\Theta'}) \hookrightarrow \mathrm{CH}(G/P_\Theta)$$

*arising from the embedding* $\Theta \subset \Theta'$ *is induced by the embedding of the respective Hasse diagrams* $\mathcal{H}_W(\Theta') \subset \mathcal{H}_W(\Theta)$.

**2.5 Corollary.** *Let $B$ be a Borel subgroup of $G$ and $P$ its (standard) parabolic subgroup. Then $\mathrm{CH}(G/P)$ is a subring of $\mathrm{CH}(G/B)$. The generators of $\mathrm{CH}(G/P)$ are $[X_w]$, where $w \in {}^\Theta W \subset W$. The cycle $[X_w]$ in $\mathrm{CH}(G/P)$ has the codimension $l(w_0) - l(w)$.*

*Proof.* Apply the lemma to the case $B = P_\emptyset$ and $P = P_{\Theta'}$. □

Hence, in order to compute $\mathrm{CH}(G/P)$ it is enough to compute $\mathrm{CH}(X)$, where $X = G/B$ is the variety of complete flags. The following results provide tools to perform such computations.

**2.6** (Poincaré duality). In order to multiply two basis elements $h$ and $g$ of $\mathrm{CH}(G/P)$ such that $\deg h + \deg g = \dim G/P$ we use the following formula (see [Ko91, 1.4]):

$$[X_w] \cdot [X_{w'}] = \delta_{w, w_0 w' w_\theta} \cdot [pt].$$

**2.7** (Pieri formula). In order to multiply two basis elements of $\mathrm{CH}(G/B)$ one of which is of codimension 1 we use the following formula (see [De74, Cor. 2 of 4.4]):

$$[X_{w_0 s_\alpha}][X_w] = \sum_{\beta \in \Phi^+,\, l(ws_\beta) = l(w)-1} \langle \beta^\vee, \bar{\omega}_\alpha \rangle [X_{ws_\beta}],$$

where the sum runs through the set of positive roots $\beta \in \Phi^+$, $s_\alpha$ denotes the simple reflection corresponding to $\alpha$ and $\bar{\omega}_\alpha$ is the fundamental weight corresponding to $\alpha$. Here $[X_{w_0 s_\alpha}]$ is the element of codimension 1.

**2.8** (Giambelli formula). Let $\mathrm{P} = \mathrm{P}(\Phi)$ be the weight space. We denote as $\bar{\omega}_1, \ldots \bar{\omega}_l$ the basis of P consisting of fundamental weights. The symmetric algebra $S^*(\mathrm{P})$ is isomorphic to $\mathbb{Z}[\bar{\omega}_1, \ldots \bar{\omega}_l]$. The Weyl group $W$ acts on P, hence, on $S^*(\mathrm{P})$. Namely, for a simple root $\alpha_i$,

$$w_{\alpha_i}(\bar{\omega}_j) = \begin{cases} \bar{\omega}_i - \alpha_i, & i = j, \\ \bar{\omega}_j, & \text{otherwise.} \end{cases}$$

We define a linear map $c\colon S^*(\mathrm{P}) \to \mathrm{CH}^*(G/B)$ as follows. For a homogeneous $u \in \mathbb{Z}[\bar{\omega}_1, \ldots, \bar{\omega}_l]$

$$c(u) = \sum_{w \in W,\, l(w) = \deg(u)} \Delta_w(u)[X_{w_0 w}],$$

where for $w = w_{\alpha_1} \ldots w_{\alpha_k}$ we denote by $\Delta_w$ the composition of derivations $\Delta_{\alpha_1} \circ \ldots \circ \Delta_{\alpha_k}$ and the derivation $\Delta_{\alpha_i}: S^*(\mathrm{P}) \to S^{*-1}(\mathrm{P})$ is defined by $\Delta_{\alpha_i}(u) = \frac{u - w_{\alpha_i}(u)}{\alpha_i}$. Then (see [Hi82a, ch. IV, 2.4])

$$[X_w] = c(\Delta_{w^{-1}}(\frac{d}{|W|})),$$

where $d$ is the product of all positive roots in $S^*(\mathrm{P})$. In other words, the element $\Delta_{w^{-1}}(\frac{d}{|W|}) \in c^{-1}([X_w])$.

Hence, in order to multiply two basis elements $h, g \in \mathrm{CH}(X)$ one may take their preimages under the map $c$ and multiply them in $S^*(\mathrm{P}) \otimes_{\mathbb{Z}} \mathbb{Q} = \mathbb{Q}[\bar{\omega}_1, \ldots \bar{\omega}_l]$, finally applying $c$ to the product.

## 3 Algorithms and implementation

The progress the world of computing is experiencing (in both software and hardware) has been so fast that we shall not be surprised if in a dozen of years computers will be able to compute with the Giambelli formula for $E_8$ as it was written down in the previous section. However, currently it is a hopeless task. In this section we show how we have nevertheless implemented the formulae described above. In some cases (namely, for $E_6$, $E_7$ and $E_8$) this required additional insight, both in mathematics and programming. The mathematical part of these insights is described here.

**3.1.** Pieri formula implementation is rather straightforward: it is quite feasible for a modern computer to search through all possible $\beta$'s (see 2.7), reduce each word and select the ones that are longer than the original even after reduction.

After implementing the basic Pieri formula — multiplication of each element by $[X_{w_0 s_\alpha}]$ — one would naturally try to exploit all information concealed in the Pieri formula and thus develop a substantial part of the multiplication table by simply taking powers of $[X_{w_0 s_\alpha}]$. At this point it also becomes clear where exactly one has to apply the Giambelli formula: in some codimension (e.g. for $X = E_8/P_8$ in codimensions 6 and 10) powers of the generator of $\mathrm{CH}^1(X)$ do not suffice.

**3.2.** The Giambelli formula as described in 2.8 is quite feasible for groups of relatively small rank. However, for $E_6$, $E_7$, and especially $E_8$ (for large classical groups also, but we are primarily concerned with covering all exceptional ones) it becomes too hard computationally. To illustrate this, we note that the formula non-trivially operates with the product of linear polynomials corresponding to all positive roots of a root system, and for $E_8$ this would be a polynomial of degree 120 in 8 variables, which is almost infeasible to write down, let alone operate with.

**3.3.** To overcome this difficulty, we introduced another method of taking preimages of the $c$ function (this is precisely the part where we needed the Giambelli formula). Note that a preimage of the $c$ function is a homogeneous polynomial of a certain (given) degree in variables $\bar{\omega}_i$ with the additional property that it is stable under the actions of the Weyl group $W_P$. In fact, to find a preimage of a certain element $x \in \mathrm{CH}^k(X)$, it is sufficient to find a polynomial $p$ such that:

- $p$ is a homogeneous polynomial of degree $k$ in $\bar{\omega}_i$;

- $\Delta_w(p) = 0$ for $w \in W \setminus {}^P W$, $l(w) = k$; alternatively one can check the condition that $p$ is $W_P$-invariant: $w_{s_i}(p) = p$ for all $i \in \Theta$; we have implemented both variants, and it turns out that the first one tends to yield more concise systems of equations, while the second one actually works faster and has allowed us to work in $E_8$;

- $c(p) = x$.

To ensure stability under $W_P$, it is sufficient to check stability under its generators. Therefore, what we did was to take a generic polynomial with undefined coefficients and solve this system of equations for these coefficients. This would be a system of linear equations.

**3.4.** Let us make an example of this technique by calculating the preimage of the generator of $\mathrm{CH}^2(\mathrm{F}_4/P_1)$.

The generic polynomial in 4 variables of degree 2 looks like the following:

```
p=a[0,0,0,2]*w[4]^2+a[0,0,1,1]*w[3]*w[4]+a[1,0,0,1]*w[1]*w[4]+
a[0,1,1,0]*w[2]*w[3]+a[0,0,2,0]*w[3]^2+a[1,0,1,0]*w[1]*w[3]+
a[2,0,0,0]*w[1]^2+a[1,1,0,0]*w[1]*w[2]+a[0,2,0,0]*w[2]^2+
a[0,1,0,1]*w[2]*w[4]
```

The first approach ($\Delta_w(p) = 0$) yields the following linear system:

$$\begin{array}{ll}
a[0,0,1,1] + a[0,0,2,0] = 0 & a[1,0,0,1] = 0 \\
a[0,1,1,0] + a[0,2,0,0] = 0 & a[0,1,0,1] = 0 \\
a[1,1,0,0] + a[0,2,0,0] = 0 & a[1,0,1,0] = 0 \\
a[0,1,1,0] + 2a[0,0,2,0] = 0 & a[0,0,0,2] + a[0,0,1,1] = 0 \\
a[2,0,0,0] + a[1,1,0,0] - 1 = 0.
\end{array}$$

After solving it, the result has one parameter that corresponds to a non-trivial polynomial of degree 2 that is stable under $\Delta_w$. The result is as follows:

$$\begin{array}{ll}
a[2,0,0,0] = 1 + 2t, & a[1,1,0,0] = -2t, \\
a[1,0,1,0] = 0, & a[1,0,0,1] = 0, \\
a[0,2,0,0] = 2t, & a[0,1,1,0] = -2t, \\
a[0,1,0,1] = 0, & a[0,0,2,0] = t, \\
a[0,0,1,1] = -t, & a[0,0,0,2] = t.
\end{array}$$

For example, a possible solution $p = \bar{\omega}_1^2$ arises from this result.

The other, equivalent checking procedure, namely applying each of the $w_{s_i}$, $i = 2, 3, 4$ is equivalent to replacing $\bar{\omega}_i$ with $\bar{\omega}_i - \alpha_i$; in the case of $F_4$ these polynomials are (we need only the last three):

$$\bar{\omega}_1 - \alpha_1 = -\bar{\omega}_1 + \bar{\omega}_2 \qquad \bar{\omega}_2 - \alpha_2 = \bar{\omega}_1 - \bar{\omega}_2 + \bar{\omega}_3$$
$$\bar{\omega}_3 - \alpha_3 = 2\bar{\omega}_2 - \bar{\omega}_3 + \bar{\omega}_4 \qquad \bar{\omega}_4 - \alpha_4 = \bar{\omega}_3 - \bar{\omega}_4$$

After the substitutions the polynomial changes as follows:

```
w[s2](p)=a[0,0,0,2]*w[4]^2+a[0,0,1,1]*w[3]*w[4]+a[1,0,0,1]
*w[1]*w[4]-a[0,1,1,0]*w[2]*w[3]+a[0,1,1,0]*w[3]*w[1]
+a[0,1,1,0]*w[3]^2+a[0,0,2,0]*w[3]^2+a[1,0,1,0]*w[1]*w[3]
+a[2,0,0,0]*w[1]^2-a[1,1,0,0]*w[1]*w[2]+a[1,1,0,0]*w[1]^2
+a[1,1,0,0]*w[1]*w[3]+a[0,2,0,0]*w[2]^2-2*a[0,2,0,0]*w[1]*w[2]
-2*a[0,2,0,0]*w[2]*w[3]+a[0,2,0,0]*w[1]^2+2*a[0,2,0,0]*w[1]
*w[3]+a[0,2,0,0]*w[3]^2-a[0,1,0,1]*w[2]*w[4]+a[0,1,0,1]
*w[4]*w[1]+a[0,1,0,1]*w[4]*w[3]
```

```
w[s3](p)=a[0,0,0,2]*w[4]^2-a[0,0,1,1]*w[3]*w[4]+2*a[0,0,1,1]
*w[4]*w[2]+a[0,0,1,1]*w[4]^2+a[1,0,0,1]*w[1]*w[4]-a[0,1,1,0]
*w[2]*w[3]+2*a[0,1,1,0]*w[2]^2+a[0,1,1,0]*w[2]*w[4]
+a[0,0,2,0]*w[3]^2-4*a[0,0,2,0]*w[2]*w[3]-2*a[0,0,2,0]*w[3]
*w[4]+4*a[0,0,2,0]*w[2]^2+4*a[0,0,2,0]*w[2]*w[4]+a[0,0,2,0]
*w[4]^2-a[1,0,1,0]*w[1]*w[3]+2*a[1,0,1,0]*w[1]*w[2]
+a[1,0,1,0]*w[1]*w[4]+a[2,0,0,0]*w[1]^2+a[1,1,0,0]*w[1]*w[2]
+a[0,2,0,0]*w[2]^2+a[0,1,0,1]*w[2]*w[4]
```

```
w[s4](p)=4*a[0,0,0,2]*w[4]^2+8*a[0,0,0,2]*w[2]*w[4]
-8*a[0,0,0,2]*w[3]*w[4]+4*a[0,0,0,2]*w[2]^2-8*a[0,0,0,2]*w[2]
*w[3]+4*a[0,0,0,2]*w[3]^2+2*a[0,0,1,1]*w[3]*w[4]+2*a[0,0,1,1]
*w[3]*w[2]-2*a[0,0,1,1]*w[3]^2+2*a[1,0,0,1]*w[1]*w[4]+2
*a[1,0,0,1]*w[1]*w[2]-2*a[1,0,0,1]*w[1]*w[3]+a[0,1,1,0]*w[2]
*w[3]+a[0,0,2,0]*w[3]^2+a[1,0,1,0]*w[1]*w[3]+a[2,0,0,0]
*w[1]^2+a[1,1,0,0]*w[1]*w[2]+a[0,2,0,0]*w[2]^2+2*a[0,1,0,1]
*w[2]*w[4]+2*a[0,1,0,1]*w[2]^2-2*a[0,1,0,1]*w[2]*w[3]
```

Now we identify the corresponding coefficients of these polynomials and get a linear system (which even in this small case is too large to display here); after solving this system, we get the same result as above.

**3.5.** After we have a polynomial $p$ from the preimage of a given element of the Chow group, we are still not quite finished. The polynomial may still be too large, especially for large codimensions, and it may still be too computationally intensive to work with it. Therefore, it is desirable to obtain a *minimal* polynomial from the preimage, or at least try to reduce the existing polynomial somehow. We have done this with the help of the Groebner bases.

Recall that $\mathrm{CH}(G/P) \otimes_{\mathbb{Z}} \mathbb{Q} \simeq \mathbb{Q}[P]^{W_P}/\mathbb{Q}[P]^{W}_+$, where $\mathbb{Q}[P]^{W_P}$ is the ring of $W_P$-invariant polynomials of the weight lattice $\mathbb{Q}[P]$ and $\mathbb{Q}[P]^{W}_+$ is the ideal of $\mathbb{Q}[P]^{W_P}$ generated by $W$-invariants without constant term (see [IM05, Sections 6]). The bases of the rings of polynomials that are invariant under different Weyl groups have already been computed (see [Meh88]). For each group we have pre-calculated the Groebner basis of the ideal of $\mathbb{Q}[P]^{W_B} = \mathbb{Q}[P]$ generated by $W$-invariants, and it is now easy to reduce the obtained polynomial $p$ by this basis. Finally, this result we accept as the preimage of the element of our Chow ring. Taking the $c$ function, not $c^{-1}$, is easier and is feasible by itself, because it does not involve the product of all roots. Apart from this, we use the Leibniz rule [Hi82a, Ch. IV, Lemma 1.1(e)] to take the $c$ function of a product of two polynomials.

**3.6.** In practice, one only needs to take very few preimages to calculate the entire Chow ring. This is so because of the Pieri formula and Poincaré duality. For example, for $X = E_8/P_8$ it suffices to find preimages only in codimensions 6 and 10.

# 4 ChowMaple: a brief outline

We refer the reader to the Web page [Chow06], where one can download the package and read a complete manual on its functionality. The package is based on John Stembridge's `coxeter` and `weyl` packages [St04]. We show here a table of `ChowMaple`'s basic functionality. The left column contains names of different functions, while the right column provides explanations. The root system type and parabolics are initialized as global variables.

| | |
|---|---|
| `init_all()` | initializes basic structures: the root system, the Weyl group, and the parabolic subgroup |
| `pieri_part()` | creates the Pieri table (calculates everything that is possible to calculate using Pieri formulae |
| `cfuncP(p)` | applies the $c$ function to a given polynomial |
| `apply_delta(w,p)` | applies the $\Delta_w$ function for a given $w$ to a given polynomial |
| `get_polynomial(c)` | finds a preimage of a given Schubert variety |
| `multchow(c1,c2)` | multiplies the two given Schubert varieties (elements of the Chow ring); implements the algorithm described in Section 3 |

# 5 Multiplicative structure of exceptional projective homogeneous varieties

In the present section we list formulae which determine the multiplicative structure of $\mathrm{CH}(E_7/P_7)$ and $\mathrm{CH}(E_8/P_8)$. The case of $F_4/P_1$ and $F_4/P_4$ was already treated in [NSZ05], where it was a substantial step on the way to the main result. Iliev and Manivel in [IM05] treated the case of $E_6/P_1$. We here finish the simplest (and, therefore, most interesting) cases, considering $E_7/P_7$ and $E_8/P_8$.

**5.1.** The basic combinatorial structure that we use is the so-called *Pieri graph*. Its nodes represent generators of Chow groups (we draw the graph in such a way that generators of the same $CH^i$ are vertically aligned, and codimension grows from right to left). As we have already noted, the generators correspond to vertices on the corresponding Hasse diagram, and the Pieri graph is thus very similar to it. The only difference is that we also, besides the labels corresponding to fundamental roots, assign weights to the edges (denoted by multiple edges on the graphs) which correspond to the coefficients in the Pieri formula. This means that if a vertex $w$ is connected with its left neighbor $w'$ with a $k$-fold edge, then $w'$ occurs in the product $wh$ with coefficient $k$, where $h$ is the canonical generator of the Chow group of codimension 1 (in all cases we are interested in there is only one such generator). Examples of Pieri graphs with multiple edges are provided in [NSZ05, Section 4].

**5.2.** Note that the case of $E_6/P_1$ was already considered by Iliev and Manivel in [IM05]. They used completely different, purely geometric techniques, making heavy use of the geometry of the Cayley plane which happens to be an $E_6/P_1$ variety. It was a great relief for us when the results produced by our package for this case completely matched the results of [IM05]. We use only combinatorics explained above, and our Maple package produces all formulae for this case in a few seconds. We refer the reader to [IM05] for the formulae themselves.

**5.3.** Consider the case $E_7/P_7$. The Pieri graph for $E_7/P_7$ looks like the following:

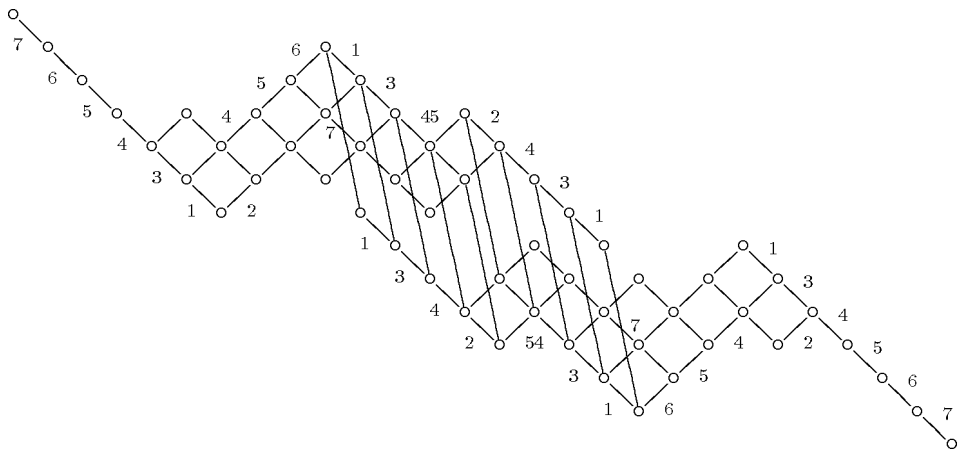

Since this case corresponds to a microweight representation, all coefficients are equal to 1, and the Pieri graph simply coincides with the respective Hasse diagram and also with the weight diagram $(E_7, \bar\omega_7)$. We note in passing that the weight diagram for $E_7/P_7$ includes two weight diagrams for $E_6/P_6$ extended with one more vertex each; for a detailed discussion of weight diagrams, their meanings and uses see [PSV98, Va00] and references therein.

We now proceed to describing the multiplicative structure of the Chow ring of $E_7/P_7$; namely, we list the generators and show how these generators are multiplied (in the $g_{i,j}$ notation translated to Weyl group words in the Appendix).

**5.4 Proposition.** *The $\mathbb{Q}$-algebra $\mathrm{CH}(E_7/P_7) \otimes_{\mathbb{Z}} \mathbb{Q}$ is generated (as an algebra) by $g_{1,1}$, $g_{5,1}$, and $g_{9,1}$. The following formulae hold:*

$$(g_{5,1})^2 = 2g_{10,1} + 2g_{10,2} \qquad (g_{5,1})^3 = 8g_{15,1} + 22g_{15,2} + 6g_{15,3}$$
$$g_{5,1}g_{9,1} = 2g_{14,1} + 3g_{14,2} + 2g_{14,3} \qquad (g_{9,1})^2 = 2g_{18,1} + 4g_{18,2} + 2g_{18,3}$$
$$(g_{5,1})^4 = 64g_{20,1} + 120g_{20,2} \qquad (g_{5,1})^5 = 184g_{25,1}$$
$$(g_{9,1})^3 = 2g_{27,1} \qquad g_{5,1}(g_{9,1})^2 = 18g_{23,1}$$
$$(g_{5,1})^2 g_{9,1} = 18g_{19,2} + 20g_{19,1} \qquad (g_{5,1})^3 g_{9,1} = 58g_{24,1}$$

and $g_{1,1}u = \sum_{v \leftarrow u} v$, where the sum runs through all the edges going from $u$ one step to the left in the Pieri graph 5.3.

**5.5.** First, let us describe the Pieri graph for $E_8/P_8$. To do this, observe that all coefficients in Pieri formulae up to the middle codimension 28 are 1's and in codimension 28 they look as follows:

$$g_{28,1}g_{1,1} = g_{29,3} + 2g_{29,1} \qquad g_{28,2}g_{1,1} = g_{29,4} + 2g_{29,2}$$
$$g_{28,3}g_{1,1} = g_{29,4} + g_{29,1} + 2g_{29,3} \qquad g_{28,4}g_{1,1} = g_{29,5} + g_{29,2} + g_{29,3} + 2g_{29,4}$$
$$g_{28,5}g_{1,1} = g_{29,4} + g_{29,6} + 2g_{29,5} \qquad g_{28,6}g_{1,1} = g_{29,5} + g_{29,7} + 2g_{29,6}$$
$$g_{28,7}g_{1,1} = g_{29,6} + g_{29,8} + 2g_{29,7} \qquad g_{28,8}g_{1,1} = g_{29,7} + 2g_{29,8}$$

In codimensions up to 28 the Pieri graph coincides with the Hasse diagram for $E_8/P_8$ and also with the weight diagram for $(E_8, \bar{\omega}_8)$ (the diagram is seldom completely drawn in literature; we recommend [PSV98, Fig. 24]) Moreover, the Pieri graph is symmetric. All above said determines the Pieri graph completely.

The multiplicative structure of the Chow ring in this case (again in the $g_{i,j}$ notation translated to Weyl group words in the Appendix) can be summarized in the following proposition.

**5.6 Proposition.** *The $\mathbb{Q}$-algebra $\mathrm{CH}(E_8/P_8) \otimes_{\mathbb{Z}} \mathbb{Q}$ is generated (as an algebra) by $g_{1,1}$, $g_{6,1}$, and $g_{10,1}$. The following formulae hold:*

$$\begin{aligned}
(g_{6,1})^2 &= 2g_{12,2} + 4g_{12,3} + 2g_{12,4} + g_{12,1}, \\
(g_{6,1})^3 &= 6g_{18,4} + 34g_{18,5} + 58g_{18,6} + 85g_{18,2} + 111g_{18,3} + 25g_{18,1}, \\
(g_{6,1})^4 &= 432g_{24,7} + 2256g_{24,5} + 1668g_{24,6} + 5600g_{24,4} + \\
&\quad + 3957g_{24,1} + 2048g_{24,3} + 2888g_{24,2}, \\
(g_{6,1})^5 &= 21260g_{30,7} + 88385g_{30,6} + 230349g_{30,5} + 372664g_{30,4} + \\
&\quad + 308080g_{30,3} + 331300g_{30,2} + 150705g_{30,1}, \\
(g_{6,1})^6 &= 13289373g_{36,1} + 13168365g_{36,2} + 9548378g_{36,3} + \\
&\quad + 27940475g_{36,4} + 13629290g_{36,5} + 21251025g_{36,6}, \\
(g_{6,1})^7 &= 349931688g_{42,1} + 1122626055g_{42,2} + 820427553g_{42,3} + \\
&\quad + 1182401805g_{42,4}, \\
(g_{6,1})^8 &= 20550124104g_{48,1} + 19656865560g_{48,2}, \\
(g_{6,1})^9 &= 60757113768g_{54,1}, \\
\\
(g_{10,1})^2 &= 6g_{20,5} + 7g_{20,2} + 6g_{20,6} + 17g_{20,3} + 11g_{20,4} + 13g_{20,1}, \\
(g_{10,1})^3 &= 245g_{30,7} + 995g_{30,6} + 2541g_{30,5} + 4051g_{30,4} + 3325g_{30,3} + \\
&\quad + 3565g_{30,2} + 1605g_{30,1}, \\
(g_{10,1})^4 &= 358330g_{40,3} + 232900g_{40,2} + 440170g_{40,1} + 686440g_{40,5} + \\
&\quad + 1046280g_{40,4}, \\
(g_{10,1})^5 &= 24311090g_{50,2} + 8223140g_{50,1}, \\
\\
(g_{6,1})(g_{10,1}) &= 2g_{16,3} + 7g_{16,4} + 4g_{16,1} + 2g_{16,5} + 4g_{16,2}, \\
(g_{6,1})(g_{10,1})^2 &= 18g_{26,7} + 188g_{26,6} + 812g_{26,5} + 1001g_{26,3} + 727g_{26,4} + \\
&\quad + 557g_{26,2} + 280g_{26,1}, \\
(g_{6,1})(g_{10,1})^3 &= 231315g_{36,6} + 148460g_{36,5} + 143817g_{36,1} + 302945g_{36,4} + \\
&\quad + 103652g_{36,3} + 142395g_{36,2}, \\
(g_{6,1})(g_{10,1})^4 &= 24311090g_{46,1} + 4632270g_{46,3} + 16827980g_{46,2}, \\
(g_{6,1})(g_{10,1})^5 &= 24311090g_{56,1}, \\
(g_{6,1})^2(g_{10,1}) &= 58g_{22,7} + 251g_{22,4} + 201g_{22,2} + 20g_{22,6} + 137g_{22,5} + \\
&\quad + 306g_{22,3} + 85g_{22,1}, \\
(g_{6,1})^2(g_{10,1})^2 &= 41165g_{32,1} + 43609g_{32,2} + 48323g_{32,4} + 47489g_{32,5} + \\
&\quad + 74550g_{32,3} + 34339g_{32,6} + 16139g_{32,7}, \\
(g_{6,1})^2(g_{10,1})^3 &= 3790422g_{42,1} + 12817545g_{42,4} + 8894817g_{42,3} + \\
&\quad + 12176265g_{42,2}, \\
(g_{6,1})^2(g_{10,1})^4 &= 193499640g_{52,1},
\end{aligned}$$

$$\begin{aligned}
(g_{6,1})^3(g_{10,1}) &= 971g_{28,1} + 19308g_{28,4} + 7270g_{28,3} + 11984g_{28,5} + \\
&\quad + 4382g_{28,2} + 4038g_{28,6} + 58g_{28,8} + 738g_{28,7}, \\
(g_{6,1})^3(g_{10,1})^2 &= 4365020g_{38,3} + 974954g_{38,1} + 2226414g_{38,2} + \\
&\quad + 5053390g_{38,5} + 3838670g_{38,6} + 1210006g_{38,4}, \\
(g_{6,1})^3(g_{10,1})^3 &= 213093030g_{48,2} + 222792846g_{48,1}, \\
(g_{6,1})^4(g_{10,1}) &= 593807g_{34,4} + 1138388g_{34,5} + 1436181g_{34,3} + \\
&\quad + 797276g_{34,1} + 833852g_{34,2} + 448604g_{34,7} + \\
&\quad + 451492g_{34,6}, \\
(g_{6,1})^4(g_{10,1})^2 &= 99724478g_{44,4} + 243616314g_{44,3} + 71721616g_{44,2} + \\
&\quad + 12906362g_{44,1}, \\
(g_{6,1})^4(g_{10,1})^3 &= 658678722g_{54,1}, \\
(g_{6,1})^5(g_{10,1}) &= 40626835g_{40,1} + 21505960g_{40,2} + 33002530g_{40,3} + \\
&\quad + 96448185g_{40,4} + 63303340g_{40,5}, \\
(g_{6,1})^5(g_{10,1})^2 &= 2242207070g_{50,2} + 758403200g_{50,1}, \\
(g_{6,1})^6(g_{10,1}) &= 2242207070g_{46,1} + 1552248110g_{46,2} + 427403580g_{46,3}, \\
(g_{6,1})^6(g_{10,1})^2 &= 2242207070g_{56,1}, \\
(g_{6,1})^7(g_{10,1}) &= 17847431370g_{52,1},
\end{aligned}$$

and $g_{1,1}u = \sum_{v \leftarrow u} v$, where the sum runs through all the edges going from $u$ one step to the left in the Pieri graph for $E_8/P_8$.

**Acknowledgements.** Our joint work has been made possible by the hospitality of the Bielefeld University and personally Anthony Bak, to whom we express our sincere gratitude.

# Appendix

Translating the $g$ notation for $E_7/P_7$ into words to be read on the weight diagram:

$g_{0,1} = [7,6,5,4,3,2,4,5,6,7,1,3,4,5,6,2,4,5,3,4,1,3,2,4,5,6,7]$,
$g_{1,1} = [7,6,5,4,3,2,4,5,6,7,1,3,4,5,6,2,4,5,3,4,1,3,2,4,5,6]$,
$g_{2,1} = [7,6,5,4,3,2,4,5,6,7,1,3,4,5,6,2,4,5,3,4,1,3,2,4,5]$,
$g_{3,1} = [7,6,5,4,3,2,4,5,6,7,1,3,4,5,6,2,4,5,3,4,1,3,2,4]$,
$g_{4,1} = [7,6,5,4,3,2,4,5,6,7,1,3,4,5,6,2,4,5,3,4,1,3,2]$,
$g_{5,1} = [7,6,5,4,3,2,4,5,6,7,1,3,4,5,6,2,4,5,3,4,1,3]$,
$g_{5,2} = [7,6,5,4,3,2,4,5,6,7,1,3,4,5,6,2,4,5,3,4,2,1]$,
$g_{6,1} = [7,6,5,4,3,2,4,5,6,7,1,3,4,5,6,2,4,5,3,4,2]$,
$g_{6,2} = [7,6,5,4,3,2,4,5,6,7,1,3,4,5,6,2,4,5,3,4,1]$,
$g_{7,1} = [7,6,5,4,3,2,4,5,6,7,1,3,4,5,6,2,4,5,3,4]$,
$g_{7,2} = [7,6,5,4,3,2,4,5,6,7,1,3,4,5,6,2,4,5,3,1]$,
$g_{8,1} = [7,6,5,4,3,2,4,5,6,7,1,3,4,5,6,2,4,5,3]$,
$g_{8,2} = [7,6,5,4,3,2,4,5,6,7,1,3,4,5,6,2,4,3,1]$,
$g_{9,1} = [7,6,5,4,3,2,4,5,6,7,1,3,4,5,6,2,4,5]$,
$g_{9,2} = [7,6,5,4,3,2,4,5,6,7,1,3,4,5,6,2,4,3]$,
$g_{9,3} = [7,6,5,4,3,2,4,5,6,7,1,3,4,5,2,4,3,1]$,
$g_{10,1} = [7,6,5,4,3,2,4,5,6,7,1,3,4,5,6,2,4]$,
$g_{10,2} = [7,6,5,4,3,2,4,5,6,7,1,3,4,5,2,4,3]$,
$g_{10,3} = [7,6,5,4,3,2,4,5,6,1,3,4,5,2,4,3,1]$,
$g_{11,1} = [7,6,5,4,3,2,4,5,6,7,1,3,4,5,6,2]$,
$g_{11,2} = [7,6,5,4,3,2,4,5,6,7,1,3,4,5,2,4]$,
$g_{11,3} = [7,6,5,4,3,2,4,5,6,1,3,4,5,2,4,3]$,
$g_{12,1} = [7,6,5,4,3,2,4,5,6,7,1,3,4,5,6]$,
$g_{12,2} = [7,6,5,4,3,2,4,5,6,7,1,3,4,5,2]$,
$g_{12,3} = [7,6,5,4,3,2,4,5,6,1,3,4,5,2,4]$,
$g_{13,1} = [7,6,5,4,3,2,4,5,6,7,1,3,4,5]$,
$g_{13,2} = [7,6,5,4,3,2,4,5,6,7,1,3,4,2]$,
$g_{13,3} = [7,6,5,4,3,2,4,5,6,1,3,4,5,2]$,
$g_{14,1} = [7,6,5,4,3,2,4,5,6,1,3,4,2]$,
$g_{14,2} = [7,6,5,4,3,2,4,5,6,1,3,4,5]$,
$g_{14,3} = [7,6,5,4,3,2,4,5,6,7,1,3,4]$,
$g_{15,1} = [7,6,5,4,3,2,4,5,1,3,4,2]$,
$g_{15,2} = [7,6,5,4,3,2,4,5,6,1,3,4]$,
$g_{15,3} = [7,6,5,4,3,2,4,5,6,7,1,3]$
$g_{16,1} = [7,6,5,4,3,2,4,5,1,3,4]$,
$g_{16,2} = [7,6,5,4,3,2,4,5,6,1,3]$,
$g_{16,3} = [7,6,5,4,3,2,4,5,6,7,1]$,
$g_{17,1} = [7,6,5,4,3,2,4,5,1,3]$,
$g_{17,2} = [7,6,5,4,3,2,4,5,6,1]$,
$g_{17,3} = [7,6,5,4,3,2,4,5,6,7]$,
$g_{18,1} = [7,6,5,4,3,2,4,1,3]$,
$g_{18,2} = [7,6,5,4,3,2,4,5,1]$,
$g_{18,3} = [7,6,5,4,3,2,4,5,6]$
$g_{19,1} = [7,6,5,4,3,2,4,1]$,
$g_{19,2} = [7,6,5,4,3,2,4,5]$,

$$g_{20,1} = [7, 6, 5, 4, 3, 2, 1],$$
$$g_{20,2} = [7, 6, 5, 4, 3, 2, 4],$$
$$g_{21,1} = [7, 6, 5, 4, 3, 1],$$
$$g_{21,2} = [7, 6, 5, 4, 3, 2],$$
$$g_{22,1} = [7, 6, 5, 4, 2],$$
$$g_{22,2} = [7, 6, 5, 4, 3],$$
$$g_{23,1} = [7, 6, 5, 4],$$
$$g_{24,1} = [7, 6, 5],$$
$$g_{25,1} = [7, 6],$$
$$g_{26,1} = [7],$$
$$g_{27,1} = [].$$

And the same for $\mathrm{E}_8/P_8$:

$g_{0,1} = [8, 7, 6, 5, 4, 3, 2, 4, 5, 6, 7, 1, 3, 4, 5, 6, 2, 4, 5, 3, 4, 1, 3, 2, 4, 5, 6, 7, 8, 7, 6, 5, 4, 3, 2, 4, 5, 6, 7, 1, 3, 4, 5, 6, 2, 4, 5, 3, 4, 1, 3, 2, 4, 5, 6, 7, 8]$

$g_{1,1} = [8, 7, 6, 5, 4, 3, 2, 4, 5, 6, 7, 1, 3, 4, 5, 6, 2, 4, 5, 3, 4, 1, 3, 2, 4, 5, 6, 7, 8, 7, 6, 5, 4, 3, 2, 4, 5, 6, 7, 1, 3, 4, 5, 6, 2, 4, 5, 3, 4, 1, 3, 2, 4, 5, 6, 7]$

$g_{2,1} = [8, 7, 6, 5, 4, 3, 2, 4, 5, 6, 7, 1, 3, 4, 5, 6, 2, 4, 5, 3, 4, 1, 3, 2, 4, 5, 6, 7, 8, 7, 6, 5, 4, 3, 2, 4, 5, 6, 7, 1, 3, 4, 5, 6, 2, 4, 5, 3, 4, 1, 3, 2, 4, 5, 6]$

$g_{3,1} = [8, 7, 6, 5, 4, 3, 2, 4, 5, 6, 7, 1, 3, 4, 5, 6, 2, 4, 5, 3, 4, 1, 3, 2, 4, 5, 6, 7, 8, 7, 6, 5, 4, 3, 2, 4, 5, 6, 7, 1, 3, 4, 5, 6, 2, 4, 5, 3, 4, 1, 3, 2, 4, 5]$

$g_{4,1} = [8, 7, 6, 5, 4, 3, 2, 4, 5, 6, 7, 1, 3, 4, 5, 6, 2, 4, 5, 3, 4, 1, 3, 2, 4, 5, 6, 7, 8, 7, 6, 5, 4, 3, 2, 4, 5, 6, 7, 1, 3, 4, 5, 6, 2, 4, 5, 3, 4, 1, 3, 2, 4]$

$g_{5,1} = [8, 7, 6, 5, 4, 3, 2, 4, 5, 6, 7, 1, 3, 4, 5, 6, 2, 4, 5, 3, 4, 1, 3, 2, 4, 5, 6, 7, 8, 7, 6, 5, 4, 3, 2, 4, 5, 6, 7, 1, 3, 4, 5, 6, 2, 4, 5, 3, 4, 1, 3, 2]$

$g_{6,1} = [8, 7, 6, 5, 4, 3, 2, 4, 5, 6, 7, 1, 3, 4, 5, 6, 2, 4, 5, 3, 4, 1, 3, 2, 4, 5, 6, 7, 8, 7, 6, 5, 4, 3, 2, 4, 5, 6, 7, 1, 3, 4, 5, 6, 2, 4, 5, 3, 4, 1, 3]$

$g_{6,2} = [8, 7, 6, 5, 4, 3, 2, 4, 5, 6, 7, 1, 3, 4, 5, 6, 2, 4, 5, 3, 4, 1, 3, 2, 4, 5, 6, 7, 8, 7, 6, 5, 4, 3, 2, 4, 5, 6, 7, 1, 3, 4, 5, 6, 2, 4, 5, 3, 4, 2, 1]$

$g_{7,1} = [8, 7, 6, 5, 4, 3, 2, 4, 5, 6, 7, 1, 3, 4, 5, 6, 2, 4, 5, 3, 4, 1, 3, 2, 4, 5, 6, 7, 8, 7, 6, 5, 4, 3, 2, 4, 5, 6, 7, 1, 3, 4, 5, 6, 2, 4, 5, 3, 4, 2]$

$g_{7,2} = [8, 7, 6, 5, 4, 3, 2, 4, 5, 6, 7, 1, 3, 4, 5, 6, 2, 4, 5, 3, 4, 1, 3, 2, 4, 5, 6, 7, 8, 7, 6, 5, 4, 3, 2, 4, 5, 6, 7, 1, 3, 4, 5, 6, 2, 4, 5, 3, 4, 1]$

$g_{8,1} = [8, 7, 6, 5, 4, 3, 2, 4, 5, 6, 7, 1, 3, 4, 5, 6, 2, 4, 5, 3, 4, 1, 3, 2, 4, 5, 6, 7, 8, 7, 6, 5, 4, 3, 2, 4, 5, 6, 7, 1, 3, 4, 5, 6, 2, 4, 5, 3, 4]$

$g_{8,2} = [8, 7, 6, 5, 4, 3, 2, 4, 5, 6, 7, 1, 3, 4, 5, 6, 2, 4, 5, 3, 4, 1, 3, 2, 4, 5, 6, 7, 8, 7, 6, 5, 4, 3, 2, 4, 5, 6, 7, 1, 3, 4, 5, 6, 2, 4, 5, 3, 1]$

$g_{9,1} = [8, 7, 6, 5, 4, 3, 2, 4, 5, 6, 7, 1, 3, 4, 5, 6, 2, 4, 5, 3, 4, 1, 3, 2, 4, 5, 6, 7, 8, 7, 6, 5, 4, 3, 2, 4, 5, 6, 7, 1, 3, 4, 5, 6, 2, 4, 5, 3]$

$g_{9,2} = [8, 7, 6, 5, 4, 3, 2, 4, 5, 6, 7, 1, 3, 4, 5, 6, 2, 4, 5, 3, 4, 1, 3, 2, 4, 5, 6, 7, 8, 7, 6, 5, 4, 3, 2, 4, 5, 6, 7, 1, 3, 4, 5, 6, 2, 4, 3, 1]$

$g_{10,1} = [8, 7, 6, 5, 4, 3, 2, 4, 5, 6, 7, 1, 3, 4, 5, 6, 2, 4, 5, 3, 4, 1, 3, 2, 4, 5, 6, 7, 8, 7, 6, 5, 4, 3, 2, 4, 5, 6, 7, 1, 3, 4, 5, 6, 2, 4, 5]$

$g_{10,2} = [8, 7, 6, 5, 4, 3, 2, 4, 5, 6, 7, 1, 3, 4, 5, 6, 2, 4, 5, 3, 4, 1, 3, 2, 4, 5, 6, 7, 8, 7, 6, 5, 4, 3, 2, 4, 5, 6, 7, 1, 3, 4, 5, 6, 2, 4, 3]$

$g_{10,3} = [8, 7, 6, 5, 4, 3, 2, 4, 5, 6, 7, 1, 3, 4, 5, 6, 2, 4, 5, 3, 4, 1, 3, 2, 4, 5, 6, 7, 8, 7, 6, 5, 4, 3, 2, 4, 5, 6, 7, 1, 3, 4, 5, 2, 4, 3, 1]$

$g_{11,1} = [8, 7, 6, 5, 4, 3, 2, 4, 5, 6, 7, 1, 3, 4, 5, 6, 2, 4, 5, 3, 4, 1, 3, 2, 4, 5, 6, 7, 8, 7, 6, 5, 4, 3, 2, 4, 5, 6, 7, 1, 3, 4, 5, 6, 2, 4]$

$g_{11,2} = [8, 7, 6, 5, 4, 3, 2, 4, 5, 6, 7, 1, 3, 4, 5, 6, 2, 4, 5, 3, 4, 1, 3, 2, 4, 5, 6, 7, 8, 7, 6, 5, 4, 3, 2, 4, 5, 6, 7, 1, 3, 4, 5, 2, 4, 3]$

$g_{11,3} = [8, 7, 6, 5, 4, 3, 2, 4, 5, 6, 7, 1, 3, 4, 5, 6, 2, 4, 5, 3, 4, 1, 3, 2, 4, 5, 6, 7, 8, 7, 6, 5, 4, 3, 2, 4, 5, 6, 1, 3, 4, 5, 2, 4, 3, 1]$

$g_{12,1} = [8, 7, 6, 5, 4, 3, 2, 4, 5, 6, 7, 8, 1, 3, 4, 5, 6, 2, 4, 5, 3, 4, 1, 3, 2, 4, 5, 6, 7, 6, 5, 4, 3, 2, 4, 5, 6, 1, 3, 4, 5, 2, 4, 3, 1]$

$g_{12,2} = [8, 7, 6, 5, 4, 3, 2, 4, 5, 6, 7, 1, 3, 4, 5, 6, 2, 4, 5, 3, 4, 1, 3, 2, 4, 5, 6, 7, 8, 7, 6, 5, 4, 3, 2, 4, 5, 6, 7, 1, 3, 4, 5, 6, 2]$

$g_{12,3} = [8, 7, 6, 5, 4, 3, 2, 4, 5, 6, 7, 1, 3, 4, 5, 6, 2, 4, 5, 3, 4, 1, 3, 2, 4, 5, 6, 7, 8, 7, 6, 5, 4, 3, 2, 4, 5, 6, 7, 1, 3, 4, 5, 2, 4]$

$g_{12,4} = [8, 7, 6, 5, 4, 3, 2, 4, 5, 6, 7, 1, 3, 4, 5, 6, 2, 4, 5, 3, 4, 1, 3, 2, 4, 5, 6, 7, 8, 7, 6, 5, 4, 3, 2, 4, 5, 6, 1, 3, 4, 5, 2, 4, 3]$

$g_{13,1} = [8, 7, 6, 5, 4, 3, 2, 4, 5, 6, 7, 8, 1, 3, 4, 5, 6, 2, 4, 5, 3, 4, 1, 3, 2, 4, 5, 6, 7, 6, 5, 4, 3, 2, 4, 5, 6, 1, 3, 4, 5, 2, 4, 3]$

$g_{13,2} = [8, 7, 6, 5, 4, 3, 2, 4, 5, 6, 7, 1, 3, 4, 5, 6, 2, 4, 5, 3, 4, 1, 3, 2, 4, 5, 6, 7, 8, 7, 6, 5, 4, 3, 2, 4, 5, 6, 7, 1, 3, 4, 5, 6]$

$g_{13,3} = [8, 7, 6, 5, 4, 3, 2, 4, 5, 6, 7, 1, 3, 4, 5, 6, 2, 4, 5, 3, 4, 1, 3, 2, 4, 5, 6, 7, 8, 7, 6, 5, 4, 3, 2, 4, 5, 6, 7, 1, 3, 4, 5, 2]$

$g_{13,4} = [8, 7, 6, 5, 4, 3, 2, 4, 5, 6, 7, 1, 3, 4, 5, 6, 2, 4, 5, 3, 4, 1, 3, 2, 4, 5, 6, 7, 8, 7, 6, 5, 4, 3, 2, 4, 5, 6, 1, 3, 4, 5, 2, 4]$

$g_{14,1} = [8, 7, 6, 5, 4, 3, 2, 4, 5, 6, 7, 8, 1, 3, 4, 5, 6, 2, 4, 5, 3, 4, 1, 3, 2, 4, 5, 6, 7, 6, 5, 4, 3, 2, 4, 5, 6, 1, 3, 4, 5, 2, 4]$

$g_{14,2} = [8, 7, 6, 5, 4, 3, 2, 4, 5, 6, 7, 1, 3, 4, 5, 6, 2, 4, 5, 3, 4, 1, 3, 2, 4, 5, 6, 7, 8, 7, 6, 5, 4, 3, 2, 4, 5, 6, 7, 1, 3, 4, 5]$

$g_{14,3} = [8, 7, 6, 5, 4, 3, 2, 4, 5, 6, 7, 1, 3, 4, 5, 6, 2, 4, 5, 3, 4, 1, 3, 2, 4, 5, 6, 7, 8, 7, 6, 5, 4, 3, 2, 4, 5, 6, 7, 1, 3, 4, 2]$

$g_{14,4} = [8, 7, 6, 5, 4, 3, 2, 4, 5, 6, 7, 1, 3, 4, 5, 6, 2, 4, 5, 3, 4, 1, 3, 2, 4, 5, 6, 7, 8, 7, 6, 5, 4, 3, 2, 4, 5, 6, 1, 3, 4, 5, 2]$

$g_{15,1} = [8, 7, 6, 5, 4, 3, 2, 4, 5, 6, 7, 8, 1, 3, 4, 5, 6, 2, 4, 5, 3, 4, 1, 3, 2, 4, 5, 6, 7, 6, 5, 4, 3, 2, 4, 5, 6, 1, 3, 4, 5, 2]$

$g_{15,2} = [8, 7, 6, 5, 4, 3, 2, 4, 5, 6, 7, 1, 3, 4, 5, 6, 2, 4, 5, 3, 4, 1, 3, 2, 4, 5, 6, 7, 8, 7, 6, 5, 4, 3, 2, 4, 5, 6, 7, 1, 3, 4]$

$g_{15,3} = [8, 7, 6, 5, 4, 3, 2, 4, 5, 6, 7, 1, 3, 4, 5, 6, 2, 4, 5, 3, 4, 1, 3, 2, 4, 5, 6, 7, 8, 7, 6, 5, 4, 3, 2, 4, 5, 6, 1, 3, 4, 5]$

$g_{15,4} = [8, 7, 6, 5, 4, 3, 2, 4, 5, 6, 7, 1, 3, 4, 5, 6, 2, 4, 5, 3, 4, 1, 3, 2, 4, 5, 6, 7, 8, 7, 6, 5, 4, 3, 2, 4, 5, 6, 1, 3, 4, 2]$

$g_{16,1} = [8, 7, 6, 5, 4, 3, 2, 4, 5, 6, 7, 8, 1, 3, 4, 5, 6, 2, 4, 5, 3, 4, 1, 3, 2, 4, 5, 6, 7, 6, 5, 4, 3, 2, 4, 5, 6, 1, 3, 4, 5]$

$g_{16,2} = [8, 7, 6, 5, 4, 3, 2, 4, 5, 6, 7, 8, 1, 3, 4, 5, 6, 2, 4, 5, 3, 4, 1, 3, 2, 4, 5, 6, 7, 6, 5, 4, 3, 2, 4, 5, 6, 1, 3, 4, 2]$

$g_{16,3} = [8, 7, 6, 5, 4, 3, 2, 4, 5, 6, 7, 1, 3, 4, 5, 6, 2, 4, 5, 3, 4, 1, 3, 2, 4, 5, 6, 7, 8, 7, 6, 5, 4, 3, 2, 4, 5, 6, 7, 1, 3]$

$g_{16,4} = [8, 7, 6, 5, 4, 3, 2, 4, 5, 6, 7, 1, 3, 4, 5, 6, 2, 4, 5, 3, 4, 1, 3, 2, 4, 5, 6, 7, 8, 7, 6, 5, 4, 3, 2, 4, 5, 6, 1, 3, 4]$

$g_{16,5} = [8, 7, 6, 5, 4, 3, 2, 4, 5, 6, 7, 1, 3, 4, 5, 6, 2, 4, 5, 3, 4, 1, 3, 2, 4, 5, 6, 7, 8, 7, 6, 5, 4, 3, 2, 4, 5, 1, 3, 4, 2]$

$g_{17,1} = [8, 7, 6, 5, 4, 3, 2, 4, 5, 6, 7, 8, 1, 3, 4, 5, 6, 2, 4, 5, 3, 4, 1, 3, 2, 4, 5, 6, 7, 6, 5, 4, 3, 2, 4, 5, 6, 1, 3, 4]$

$g_{17,2} = [8, 7, 6, 5, 4, 3, 2, 4, 5, 6, 7, 8, 1, 3, 4, 5, 6, 2, 4, 5, 3, 4, 1, 3, 2, 4, 5, 6, 7, 6, 5, 4, 3, 2, 4, 5, 1, 3, 4, 2]$

$g_{17,3} = [8, 7, 6, 5, 4, 3, 2, 4, 5, 6, 7, 1, 3, 4, 5, 6, 2, 4, 5, 3, 4, 1, 3, 2, 4, 5, 6, 7, 8, 7, 6, 5, 4, 3, 2, 4, 5, 6, 7, 1]$

$g_{17,4} = [8, 7, 6, 5, 4, 3, 2, 4, 5, 6, 7, 1, 3, 4, 5, 6, 2, 4, 5, 3, 4, 1, 3, 2, 4, 5, 6, 7, 8, 7, 6, 5, 4, 3, 2, 4, 5, 6, 1, 3]$

$g_{17,5} = [8, 7, 6, 5, 4, 3, 2, 4, 5, 6, 7, 1, 3, 4, 5, 6, 2, 4, 5, 3, 4, 1, 3, 2, 4, 5, 6, 7, 8, 7, 6, 5, 4, 3, 2, 4, 5, 1, 3, 4]$

$g_{18,1} = [8, 7, 6, 5, 4, 3, 2, 4, 5, 6, 7, 8, 1, 3, 4, 5, 6, 7, 2, 4, 5, 3, 4, 1, 3, 2, 4, 5, 6, 5, 4, 3, 2, 4, 5, 1, 3, 4, 2]$

$g_{18,2} = [8, 7, 6, 5, 4, 3, 2, 4, 5, 6, 7, 8, 1, 3, 4, 5, 6, 2, 4, 5, 3, 4, 1, 3, 2, 4, 5, 6, 7, 6, 5, 4, 3, 2, 4, 5, 6, 1, 3]$

$g_{18,3} = [8, 7, 6, 5, 4, 3, 2, 4, 5, 6, 7, 8, 1, 3, 4, 5, 6, 2, 4, 5, 3, 4, 1, 3, 2, 4, 5, 6, 7, 6, 5, 4, 3, 2, 4, 5, 1, 3, 4]$

$g_{18,4} = [8, 7, 6, 5, 4, 3, 2, 4, 5, 6, 7, 1, 3, 4, 5, 6, 2, 4, 5, 3, 4, 1, 3, 2, 4, 5, 6, 7, 8, 7, 6, 5, 4, 3, 2, 4, 5, 6, 7]$

$g_{18,5} = [8, 7, 6, 5, 4, 3, 2, 4, 5, 6, 7, 1, 3, 4, 5, 6, 2, 4, 5, 3, 4, 1, 3, 2, 4, 5, 6, 7, 8, 7, 6, 5, 4, 3, 2, 4, 5, 6, 1]$

$g_{18,6} = [8, 7, 6, 5, 4, 3, 2, 4, 5, 6, 7, 1, 3, 4, 5, 6, 2, 4, 5, 3, 4, 1, 3, 2, 4, 5, 6, 7, 8, 7, 6, 5, 4, 3, 2, 4, 5, 1, 3]$

$g_{19,1} = [8, 7, 6, 5, 4, 3, 2, 4, 5, 6, 7, 8, 1, 3, 4, 5, 6, 7, 2, 4, 5, 3, 4, 1, 3, 2, 4, 5, 6, 5, 4, 3, 2, 4, 5, 1, 3, 4]$

$g_{19,2} = [8, 7, 6, 5, 4, 3, 2, 4, 5, 6, 7, 8, 1, 3, 4, 5, 6, 2, 4, 5, 3, 4, 1, 3, 2, 4, 5, 6, 7, 6, 5, 4, 3, 2, 4, 5, 6, 1]$

$g_{19,3} = [8, 7, 6, 5, 4, 3, 2, 4, 5, 6, 7, 8, 1, 3, 4, 5, 6, 2, 4, 5, 3, 4, 1, 3, 2, 4, 5, 6, 7, 6, 5, 4, 3, 2, 4, 5, 1, 3]$

$g_{19,4} = [8, 7, 6, 5, 4, 3, 2, 4, 5, 6, 7, 1, 3, 4, 5, 6, 2, 4, 5, 3, 4, 1, 3, 2, 4, 5, 6, 7, 8, 7, 6, 5, 4, 3, 2, 4, 5, 6]$

$g_{19,5} = [8, 7, 6, 5, 4, 3, 2, 4, 5, 6, 7, 1, 3, 4, 5, 6, 2, 4, 5, 3, 4, 1, 3, 2, 4, 5, 6, 7, 8, 7, 6, 5, 4, 3, 2, 4, 5, 1]$

$g_{19,6} = [8, 7, 6, 5, 4, 3, 2, 4, 5, 6, 7, 1, 3, 4, 5, 6, 2, 4, 5, 3, 4, 1, 3, 2, 4, 5, 6, 7, 8, 7, 6, 5, 4, 3, 2, 4, 1, 3]$

$g_{20,1} = [8, 7, 6, 5, 4, 3, 2, 4, 5, 6, 7, 8, 1, 3, 4, 5, 6, 7, 2, 4, 5, 3, 4, 1, 3, 2, 4, 5, 6, 5, 4, 3, 2, 4, 5, 1, 3]$

$g_{20,2} = [8, 7, 6, 5, 4, 3, 2, 4, 5, 6, 7, 8, 1, 3, 4, 5, 6, 2, 4, 5, 3, 4, 1, 3, 2, 4, 5, 6, 7, 6, 5, 4, 3, 2, 4, 5, 6]$

$g_{20,3} = [8, 7, 6, 5, 4, 3, 2, 4, 5, 6, 7, 8, 1, 3, 4, 5, 6, 2, 4, 5, 3, 4, 1, 3, 2, 4, 5, 6, 7, 6, 5, 4, 3, 2, 4, 5, 1]$

$g_{20,4} = [8, 7, 6, 5, 4, 3, 2, 4, 5, 6, 7, 8, 1, 3, 4, 5, 6, 2, 4, 5, 3, 4, 1, 3, 2, 4, 5, 6, 7, 6, 5, 4, 3, 2, 4, 1, 3]$

$g_{20,5} = [8, 7, 6, 5, 4, 3, 2, 4, 5, 6, 7, 1, 3, 4, 5, 6, 2, 4, 5, 3, 4, 1, 3, 2, 4, 5, 6, 7, 8, 7, 6, 5, 4, 3, 2, 4, 5]$

$g_{20,6} = [8, 7, 6, 5, 4, 3, 2, 4, 5, 6, 7, 1, 3, 4, 5, 6, 2, 4, 5, 3, 4, 1, 3, 2, 4, 5, 6, 7, 8, 7, 6, 5, 4, 3, 2, 4, 1]$

$g_{21,1} = [8, 7, 6, 5, 4, 3, 2, 4, 5, 6, 7, 8, 1, 3, 4, 5, 6, 7, 2, 4, 5, 3, 4, 1, 3, 2, 4, 5, 6, 5, 4, 3, 2, 4, 5, 1]$

$g_{21,2} = [8, 7, 6, 5, 4, 3, 2, 4, 5, 6, 7, 8, 1, 3, 4, 5, 6, 7, 2, 4, 5, 3, 4, 1, 3, 2, 4, 5, 6, 5, 4, 3, 2, 4, 1, 3]$

$g_{21,3} = [8, 7, 6, 5, 4, 3, 2, 4, 5, 6, 7, 8, 1, 3, 4, 5, 6, 2, 4, 5, 3, 4, 1, 3, 2, 4, 5, 6, 7, 6, 5, 4, 3, 2, 4, 5]$

$g_{21,4} = [8, 7, 6, 5, 4, 3, 2, 4, 5, 6, 7, 8, 1, 3, 4, 5, 6, 2, 4, 5, 3, 4, 1, 3, 2, 4, 5, 6, 7, 6, 5, 4, 3, 2, 4, 1]$

$g_{21,5} = [8, 7, 6, 5, 4, 3, 2, 4, 5, 6, 7, 1, 3, 4, 5, 6, 2, 4, 5, 3, 4, 1, 3, 2, 4, 5, 6, 7, 8, 7, 6, 5, 4, 3, 2, 4]$

$g_{21,6} = [8, 7, 6, 5, 4, 3, 2, 4, 5, 6, 7, 1, 3, 4, 5, 6, 2, 4, 5, 3, 4, 1, 3, 2, 4, 5, 6, 7, 8, 7, 6, 5, 4, 3, 2, 1]$

$g_{22,1} = [8, 7, 6, 5, 4, 3, 2, 4, 5, 6, 7, 8, 1, 3, 4, 5, 6, 7, 2, 4, 5, 6, 3, 4, 1, 3, 2, 4, 5, 4, 3, 2, 4, 1, 3]$

$g_{22,2} = [8, 7, 6, 5, 4, 3, 2, 4, 5, 6, 7, 8, 1, 3, 4, 5, 6, 7, 2, 4, 5, 3, 4, 1, 3, 2, 4, 5, 6, 5, 4, 3, 2, 4, 5]$

$g_{22,3} = [8, 7, 6, 5, 4, 3, 2, 4, 5, 6, 7, 8, 1, 3, 4, 5, 6, 7, 2, 4, 5, 3, 4, 1, 3, 2, 4, 5, 6, 5, 4, 3, 2, 4, 1]$

$g_{22,4} = [8, 7, 6, 5, 4, 3, 2, 4, 5, 6, 7, 8, 1, 3, 4, 5, 6, 2, 4, 5, 3, 4, 1, 3, 2, 4, 5, 6, 7, 6, 5, 4, 3, 2, 4]$

$g_{22,5} = [8, 7, 6, 5, 4, 3, 2, 4, 5, 6, 7, 8, 1, 3, 4, 5, 6, 2, 4, 5, 3, 4, 1, 3, 2, 4, 5, 6, 7, 6, 5, 4, 3, 2, 1]$

$g_{22,6} = [8, 7, 6, 5, 4, 3, 2, 4, 5, 6, 7, 1, 3, 4, 5, 6, 2, 4, 5, 3, 4, 1, 3, 2, 4, 5, 6, 7, 8, 7, 6, 5, 4, 3, 1]$

$g_{22,7} = [8, 7, 6, 5, 4, 3, 2, 4, 5, 6, 7, 1, 3, 4, 5, 6, 2, 4, 5, 3, 4, 1, 3, 2, 4, 5, 6, 7, 8, 7, 6, 5, 4, 3, 2]$

$g_{23,1} = [8, 7, 6, 5, 4, 3, 2, 4, 5, 6, 7, 8, 1, 3, 4, 5, 6, 7, 2, 4, 5, 6, 3, 4, 1, 3, 2, 4, 5, 4, 3, 2, 4, 1]$
$g_{23,2} = [8, 7, 6, 5, 4, 3, 2, 4, 5, 6, 7, 8, 1, 3, 4, 5, 6, 7, 2, 4, 5, 3, 4, 1, 3, 2, 4, 5, 6, 5, 4, 3, 2, 4]$
$g_{23,3} = [8, 7, 6, 5, 4, 3, 2, 4, 5, 6, 7, 8, 1, 3, 4, 5, 6, 7, 2, 4, 5, 3, 4, 1, 3, 2, 4, 5, 6, 5, 4, 3, 2, 1]$
$g_{23,4} = [8, 7, 6, 5, 4, 3, 2, 4, 5, 6, 7, 8, 1, 3, 4, 5, 6, 2, 4, 5, 3, 4, 1, 3, 2, 4, 5, 6, 7, 6, 5, 4, 3, 1]$
$g_{23,5} = [8, 7, 6, 5, 4, 3, 2, 4, 5, 6, 7, 8, 1, 3, 4, 5, 6, 2, 4, 5, 3, 4, 1, 3, 2, 4, 5, 6, 7, 6, 5, 4, 3, 2]$
$g_{23,6} = [8, 7, 6, 5, 4, 3, 2, 4, 5, 6, 7, 1, 3, 4, 5, 6, 2, 4, 5, 3, 4, 1, 3, 2, 4, 5, 6, 7, 8, 7, 6, 5, 4, 3]$
$g_{23,7} = [8, 7, 6, 5, 4, 3, 2, 4, 5, 6, 7, 1, 3, 4, 5, 6, 2, 4, 5, 3, 4, 1, 3, 2, 4, 5, 6, 7, 8, 7, 6, 5, 4, 2]$
$g_{24,1} = [8, 7, 6, 5, 4, 3, 2, 4, 5, 6, 7, 8, 1, 3, 4, 5, 6, 7, 2, 4, 5, 6, 3, 4, 1, 3, 2, 4, 5, 4, 3, 2, 4]$
$g_{24,2} = [8, 7, 6, 5, 4, 3, 2, 4, 5, 6, 7, 8, 1, 3, 4, 5, 6, 7, 2, 4, 5, 6, 3, 4, 1, 3, 2, 4, 5, 4, 3, 2, 1]$
$g_{24,3} = [8, 7, 6, 5, 4, 3, 2, 4, 5, 6, 7, 8, 1, 3, 4, 5, 6, 7, 2, 4, 5, 3, 4, 1, 3, 2, 4, 5, 6, 5, 4, 3, 1]$
$g_{24,4} = [8, 7, 6, 5, 4, 3, 2, 4, 5, 6, 7, 8, 1, 3, 4, 5, 6, 7, 2, 4, 5, 3, 4, 1, 3, 2, 4, 5, 6, 5, 4, 3, 2]$
$g_{24,5} = [8, 7, 6, 5, 4, 3, 2, 4, 5, 6, 7, 8, 1, 3, 4, 5, 6, 2, 4, 5, 3, 4, 1, 3, 2, 4, 5, 6, 7, 6, 5, 4, 3]$
$g_{24,6} = [8, 7, 6, 5, 4, 3, 2, 4, 5, 6, 7, 8, 1, 3, 4, 5, 6, 2, 4, 5, 3, 4, 1, 3, 2, 4, 5, 6, 7, 6, 5, 4, 2]$
$g_{24,7} = [8, 7, 6, 5, 4, 3, 2, 4, 5, 6, 7, 1, 3, 4, 5, 6, 2, 4, 5, 3, 4, 1, 3, 2, 4, 5, 6, 7, 8, 7, 6, 5, 4]$
$g_{25,1} = [8, 7, 6, 5, 4, 3, 2, 4, 5, 6, 7, 8, 1, 3, 4, 5, 6, 7, 2, 4, 5, 6, 3, 4, 5, 1, 3, 2, 4, 3, 2, 1]$
$g_{25,2} = [8, 7, 6, 5, 4, 3, 2, 4, 5, 6, 7, 8, 1, 3, 4, 5, 6, 7, 2, 4, 5, 6, 3, 4, 1, 3, 2, 4, 5, 4, 3, 1]$
$g_{25,3} = [8, 7, 6, 5, 4, 3, 2, 4, 5, 6, 7, 8, 1, 3, 4, 5, 6, 7, 2, 4, 5, 6, 3, 4, 1, 3, 2, 4, 5, 4, 3, 2]$
$g_{25,4} = [8, 7, 6, 5, 4, 3, 2, 4, 5, 6, 7, 8, 1, 3, 4, 5, 6, 7, 2, 4, 5, 3, 4, 1, 3, 2, 4, 5, 6, 5, 4, 3]$
$g_{25,5} = [8, 7, 6, 5, 4, 3, 2, 4, 5, 6, 7, 8, 1, 3, 4, 5, 6, 7, 2, 4, 5, 3, 4, 1, 3, 2, 4, 5, 6, 5, 4, 2]$
$g_{25,6} = [8, 7, 6, 5, 4, 3, 2, 4, 5, 6, 7, 8, 1, 3, 4, 5, 6, 2, 4, 5, 3, 4, 1, 3, 2, 4, 5, 6, 7, 6, 5, 4]$
$g_{25,7} = [8, 7, 6, 5, 4, 3, 2, 4, 5, 6, 7, 1, 3, 4, 5, 6, 2, 4, 5, 3, 4, 1, 3, 2, 4, 5, 6, 7, 8, 7, 6, 5]$
$g_{26,1} = [8, 7, 6, 5, 4, 3, 2, 4, 5, 6, 7, 8, 1, 3, 4, 5, 6, 7, 2, 4, 5, 6, 3, 4, 5, 1, 3, 2, 4, 3, 1]$
$g_{26,2} = [8, 7, 6, 5, 4, 3, 2, 4, 5, 6, 7, 8, 1, 3, 4, 5, 6, 7, 2, 4, 5, 6, 3, 4, 5, 1, 3, 2, 4, 3, 2]$
$g_{26,3} = [8, 7, 6, 5, 4, 3, 2, 4, 5, 6, 7, 8, 1, 3, 4, 5, 6, 7, 2, 4, 5, 6, 3, 4, 1, 3, 2, 4, 5, 4, 3]$
$g_{26,4} = [8, 7, 6, 5, 4, 3, 2, 4, 5, 6, 7, 8, 1, 3, 4, 5, 6, 7, 2, 4, 5, 6, 3, 4, 1, 3, 2, 4, 5, 4, 2]$
$g_{26,5} = [8, 7, 6, 5, 4, 3, 2, 4, 5, 6, 7, 8, 1, 3, 4, 5, 6, 7, 2, 4, 5, 3, 4, 1, 3, 2, 4, 5, 6, 5, 4]$
$g_{26,6} = [8, 7, 6, 5, 4, 3, 2, 4, 5, 6, 7, 8, 1, 3, 4, 5, 6, 2, 4, 5, 3, 4, 1, 3, 2, 4, 5, 6, 7, 6, 5]$
$g_{26,7} = [8, 7, 6, 5, 4, 3, 2, 4, 5, 6, 7, 1, 3, 4, 5, 6, 2, 4, 5, 3, 4, 1, 3, 2, 4, 5, 6, 7, 8, 7, 6]$
$g_{27,1} = [8, 7, 6, 5, 4, 3, 2, 4, 5, 6, 7, 8, 1, 3, 4, 5, 6, 7, 2, 4, 5, 6, 3, 4, 5, 2, 4, 1, 3, 1]$
$g_{27,2} = [8, 7, 6, 5, 4, 3, 2, 4, 5, 6, 7, 8, 1, 3, 4, 5, 6, 7, 2, 4, 5, 6, 3, 4, 5, 1, 3, 2, 4, 3]$
$g_{27,3} = [8, 7, 6, 5, 4, 3, 2, 4, 5, 6, 7, 8, 1, 3, 4, 5, 6, 7, 2, 4, 5, 6, 3, 4, 5, 1, 3, 2, 4, 2]$
$g_{27,4} = [8, 7, 6, 5, 4, 3, 2, 4, 5, 6, 7, 8, 1, 3, 4, 5, 6, 7, 2, 4, 5, 6, 3, 4, 1, 3, 2, 4, 5, 4]$
$g_{27,5} = [8, 7, 6, 5, 4, 3, 2, 4, 5, 6, 7, 8, 1, 3, 4, 5, 6, 7, 2, 4, 5, 3, 4, 1, 3, 2, 4, 5, 6, 5]$
$g_{27,6} = [8, 7, 6, 5, 4, 3, 2, 4, 5, 6, 7, 8, 1, 3, 4, 5, 6, 2, 4, 5, 3, 4, 1, 3, 2, 4, 5, 6, 7, 6]$
$g_{27,7} = [8, 7, 6, 5, 4, 3, 2, 4, 5, 6, 7, 1, 3, 4, 5, 6, 2, 4, 5, 3, 4, 1, 3, 2, 4, 5, 6, 7, 8, 7]$
$g_{28,1} = [8, 7, 6, 5, 4, 3, 2, 4, 5, 6, 7, 8, 1, 3, 4, 5, 6, 7, 2, 4, 5, 6, 3, 4, 5, 2, 4, 3, 1]$
$g_{28,2} = [8, 7, 6, 5, 4, 3, 2, 4, 5, 6, 7, 8, 1, 3, 4, 5, 6, 7, 2, 4, 5, 6, 3, 4, 5, 1, 3, 4, 2]$
$g_{28,3} = [8, 7, 6, 5, 4, 3, 2, 4, 5, 6, 7, 8, 1, 3, 4, 5, 6, 7, 2, 4, 5, 6, 3, 4, 5, 2, 4, 1, 3]$
$g_{28,4} = [8, 7, 6, 5, 4, 3, 2, 4, 5, 6, 7, 8, 1, 3, 4, 5, 6, 7, 2, 4, 5, 6, 3, 4, 5, 1, 3, 2, 4]$
$g_{28,5} = [8, 7, 6, 5, 4, 3, 2, 4, 5, 6, 7, 8, 1, 3, 4, 5, 6, 7, 2, 4, 5, 6, 3, 4, 1, 3, 2, 4, 5]$
$g_{28,6} = [8, 7, 6, 5, 4, 3, 2, 4, 5, 6, 7, 8, 1, 3, 4, 5, 6, 7, 2, 4, 5, 3, 4, 1, 3, 2, 4, 5, 6]$
$g_{28,7} = [8, 7, 6, 5, 4, 3, 2, 4, 5, 6, 7, 8, 1, 3, 4, 5, 6, 2, 4, 5, 3, 4, 1, 3, 2, 4, 5, 6, 7]$
$g_{28,8} = [8, 7, 6, 5, 4, 3, 2, 4, 5, 6, 7, 1, 3, 4, 5, 6, 2, 4, 5, 3, 4, 1, 3, 2, 4, 5, 6, 7, 8]$
$g_{29,1} = [8, 7, 6, 5, 4, 3, 2, 4, 5, 6, 7, 8, 1, 3, 4, 5, 6, 7, 2, 4, 5, 6, 3, 4, 5, 2, 4, 3]$
$g_{29,2} = [8, 7, 6, 5, 4, 3, 2, 4, 5, 6, 7, 8, 1, 3, 4, 5, 6, 7, 2, 4, 5, 6, 3, 4, 5, 1, 3, 4]$
$g_{29,3} = [8, 7, 6, 5, 4, 3, 2, 4, 5, 6, 7, 8, 1, 3, 4, 5, 6, 7, 2, 4, 5, 6, 3, 4, 5, 2, 4, 1]$
$g_{29,4} = [8, 7, 6, 5, 4, 3, 2, 4, 5, 6, 7, 8, 1, 3, 4, 5, 6, 7, 2, 4, 5, 6, 3, 4, 5, 1, 3, 2]$
$g_{29,5} = [8, 7, 6, 5, 4, 3, 2, 4, 5, 6, 7, 8, 1, 3, 4, 5, 6, 7, 2, 4, 5, 6, 3, 4, 1, 3, 2, 4]$
$g_{29,6} = [8, 7, 6, 5, 4, 3, 2, 4, 5, 6, 7, 8, 1, 3, 4, 5, 6, 7, 2, 4, 5, 3, 4, 1, 3, 2, 4, 5]$
$g_{29,7} = [8, 7, 6, 5, 4, 3, 2, 4, 5, 6, 7, 8, 1, 3, 4, 5, 6, 2, 4, 5, 3, 4, 1, 3, 2, 4, 5, 6]$
$g_{29,8} = [8, 7, 6, 5, 4, 3, 2, 4, 5, 6, 7, 1, 3, 4, 5, 6, 2, 4, 5, 3, 4, 1, 3, 2, 4, 5, 6, 7]$
$g_{30,1} = [8, 7, 6, 5, 4, 3, 2, 4, 5, 6, 7, 8, 1, 3, 4, 5, 6, 7, 2, 4, 5, 6, 3, 4, 5, 2, 4]$
$g_{30,2} = [8, 7, 6, 5, 4, 3, 2, 4, 5, 6, 7, 8, 1, 3, 4, 5, 6, 7, 2, 4, 5, 6, 3, 4, 5, 2, 1]$
$g_{30,3} = [8, 7, 6, 5, 4, 3, 2, 4, 5, 6, 7, 8, 1, 3, 4, 5, 6, 7, 2, 4, 5, 6, 3, 4, 5, 1, 3]$
$g_{30,4} = [8, 7, 6, 5, 4, 3, 2, 4, 5, 6, 7, 8, 1, 3, 4, 5, 6, 7, 2, 4, 5, 6, 3, 4, 1, 3, 2]$
$g_{30,5} = [8, 7, 6, 5, 4, 3, 2, 4, 5, 6, 7, 8, 1, 3, 4, 5, 6, 7, 2, 4, 5, 3, 4, 1, 3, 2, 4]$
$g_{30,6} = [8, 7, 6, 5, 4, 3, 2, 4, 5, 6, 7, 8, 1, 3, 4, 5, 6, 2, 4, 5, 3, 4, 1, 3, 2, 4, 5]$
$g_{30,7} = [8, 7, 6, 5, 4, 3, 2, 4, 5, 6, 7, 1, 3, 4, 5, 6, 2, 4, 5, 3, 4, 1, 3, 2, 4, 5, 6]$
$g_{31,1} = [8, 7, 6, 5, 4, 3, 2, 4, 5, 6, 7, 8, 1, 3, 4, 5, 6, 7, 2, 4, 5, 6, 3, 4, 5, 2]$
$g_{31,2} = [8, 7, 6, 5, 4, 3, 2, 4, 5, 6, 7, 8, 1, 3, 4, 5, 6, 7, 2, 4, 5, 6, 3, 4, 5, 1]$
$g_{31,3} = [8, 7, 6, 5, 4, 3, 2, 4, 5, 6, 7, 8, 1, 3, 4, 5, 6, 7, 2, 4, 5, 6, 3, 4, 2, 1]$
$g_{31,4} = [8, 7, 6, 5, 4, 3, 2, 4, 5, 6, 7, 8, 1, 3, 4, 5, 6, 7, 2, 4, 5, 6, 3, 4, 1, 3]$
$g_{31,5} = [8, 7, 6, 5, 4, 3, 2, 4, 5, 6, 7, 8, 1, 3, 4, 5, 6, 7, 2, 4, 5, 3, 4, 1, 3, 2]$
$g_{31,6} = [8, 7, 6, 5, 4, 3, 2, 4, 5, 6, 7, 8, 1, 3, 4, 5, 6, 2, 4, 5, 3, 4, 1, 3, 2, 4]$
$g_{31,7} = [8, 7, 6, 5, 4, 3, 2, 4, 5, 6, 7, 1, 3, 4, 5, 6, 2, 4, 5, 3, 4, 1, 3, 2, 4, 5]$
$g_{32,1} = [8, 7, 6, 5, 4, 3, 2, 4, 5, 6, 7, 8, 1, 3, 4, 5, 6, 7, 2, 4, 5, 6, 3, 4, 5]$
$g_{32,2} = [8, 7, 6, 5, 4, 3, 2, 4, 5, 6, 7, 8, 1, 3, 4, 5, 6, 7, 2, 4, 5, 6, 3, 4, 2]$
$g_{32,3} = [8, 7, 6, 5, 4, 3, 2, 4, 5, 6, 7, 8, 1, 3, 4, 5, 6, 7, 2, 4, 5, 6, 3, 4, 1]$
$g_{32,4} = [8, 7, 6, 5, 4, 3, 2, 4, 5, 6, 7, 8, 1, 3, 4, 5, 6, 7, 2, 4, 5, 3, 4, 2, 1]$
$g_{32,5} = [8, 7, 6, 5, 4, 3, 2, 4, 5, 6, 7, 8, 1, 3, 4, 5, 6, 7, 2, 4, 5, 3, 4, 1, 3]$
$g_{32,6} = [8, 7, 6, 5, 4, 3, 2, 4, 5, 6, 7, 8, 1, 3, 4, 5, 6, 2, 4, 5, 3, 4, 1, 3, 2]$
$g_{32,7} = [8, 7, 6, 5, 4, 3, 2, 4, 5, 6, 7, 1, 3, 4, 5, 6, 2, 4, 5, 3, 4, 1, 3, 2, 4]$
$g_{33,1} = [8, 7, 6, 5, 4, 3, 2, 4, 5, 6, 7, 8, 1, 3, 4, 5, 6, 7, 2, 4, 5, 6, 3, 1]$
$g_{33,2} = [8, 7, 6, 5, 4, 3, 2, 4, 5, 6, 7, 8, 1, 3, 4, 5, 6, 7, 2, 4, 5, 6, 3, 4]$
$g_{33,3} = [8, 7, 6, 5, 4, 3, 2, 4, 5, 6, 7, 8, 1, 3, 4, 5, 6, 7, 2, 4, 5, 3, 4, 2]$
$g_{33,4} = [8, 7, 6, 5, 4, 3, 2, 4, 5, 6, 7, 8, 1, 3, 4, 5, 6, 7, 2, 4, 5, 3, 4, 1]$
$g_{33,5} = [8, 7, 6, 5, 4, 3, 2, 4, 5, 6, 7, 8, 1, 3, 4, 5, 6, 2, 4, 5, 3, 4, 2, 1]$
$g_{33,6} = [8, 7, 6, 5, 4, 3, 2, 4, 5, 6, 7, 8, 1, 3, 4, 5, 6, 2, 4, 5, 3, 4, 1, 3]$
$g_{33,7} = [8, 7, 6, 5, 4, 3, 2, 4, 5, 6, 7, 1, 3, 4, 5, 6, 2, 4, 5, 3, 4, 1, 3, 2]$

$g_{34,1} = [8, 7, 6, 5, 4, 3, 2, 4, 5, 6, 7, 8, 1, 3, 4, 5, 6, 7, 2, 4, 5, 6, 3]$
$g_{34,2} = [8, 7, 6, 5, 4, 3, 2, 4, 5, 6, 7, 8, 1, 3, 4, 5, 6, 7, 2, 4, 5, 3, 1]$
$g_{34,3} = [8, 7, 6, 5, 4, 3, 2, 4, 5, 6, 7, 8, 1, 3, 4, 5, 6, 7, 2, 4, 5, 3, 4]$
$g_{34,4} = [8, 7, 6, 5, 4, 3, 2, 4, 5, 6, 7, 8, 1, 3, 4, 5, 6, 2, 4, 5, 3, 4, 2]$
$g_{34,5} = [8, 7, 6, 5, 4, 3, 2, 4, 5, 6, 7, 8, 1, 3, 4, 5, 6, 2, 4, 5, 3, 4, 1]$
$g_{34,6} = [8, 7, 6, 5, 4, 3, 2, 4, 5, 6, 7, 1, 3, 4, 5, 6, 2, 4, 5, 3, 4, 2, 1]$
$g_{34,7} = [8, 7, 6, 5, 4, 3, 2, 4, 5, 6, 7, 1, 3, 4, 5, 6, 2, 4, 5, 3, 4, 1, 3]$
$g_{35,1} = [8, 7, 6, 5, 4, 3, 2, 4, 5, 6, 7, 8, 1, 3, 4, 5, 6, 7, 2, 4, 5, 6]$
$g_{35,2} = [8, 7, 6, 5, 4, 3, 2, 4, 5, 6, 7, 8, 1, 3, 4, 5, 6, 7, 2, 4, 3, 1]$
$g_{35,3} = [8, 7, 6, 5, 4, 3, 2, 4, 5, 6, 7, 8, 1, 3, 4, 5, 6, 7, 2, 4, 5, 3]$
$g_{35,4} = [8, 7, 6, 5, 4, 3, 2, 4, 5, 6, 7, 8, 1, 3, 4, 5, 6, 2, 4, 5, 3, 1]$
$g_{35,5} = [8, 7, 6, 5, 4, 3, 2, 4, 5, 6, 7, 8, 1, 3, 4, 5, 6, 2, 4, 5, 3, 4]$
$g_{35,6} = [8, 7, 6, 5, 4, 3, 2, 4, 5, 6, 7, 1, 3, 4, 5, 6, 2, 4, 5, 3, 4, 2]$
$g_{35,7} = [8, 7, 6, 5, 4, 3, 2, 4, 5, 6, 7, 1, 3, 4, 5, 6, 2, 4, 5, 3, 4, 1]$
$g_{36,1} = [8, 7, 6, 5, 4, 3, 2, 4, 5, 6, 7, 8, 1, 3, 4, 5, 6, 7, 2, 4, 3]$
$g_{36,2} = [8, 7, 6, 5, 4, 3, 2, 4, 5, 6, 7, 8, 1, 3, 4, 5, 6, 7, 2, 4, 5]$
$g_{36,3} = [8, 7, 6, 5, 4, 3, 2, 4, 5, 6, 7, 8, 1, 3, 4, 5, 6, 2, 4, 3, 1]$
$g_{36,4} = [8, 7, 6, 5, 4, 3, 2, 4, 5, 6, 7, 8, 1, 3, 4, 5, 6, 2, 4, 5, 3]$
$g_{36,5} = [8, 7, 6, 5, 4, 3, 2, 4, 5, 6, 7, 1, 3, 4, 5, 6, 2, 4, 5, 3, 1]$
$g_{36,6} = [8, 7, 6, 5, 4, 3, 2, 4, 5, 6, 7, 1, 3, 4, 5, 6, 2, 4, 5, 3, 4]$
$g_{37,1} = [8, 7, 6, 5, 4, 3, 2, 4, 5, 6, 7, 8, 1, 3, 4, 5, 6, 7, 2, 4]$
$g_{37,2} = [8, 7, 6, 5, 4, 3, 2, 4, 5, 6, 7, 8, 1, 3, 4, 5, 2, 4, 3, 1]$
$g_{37,3} = [8, 7, 6, 5, 4, 3, 2, 4, 5, 6, 7, 8, 1, 3, 4, 5, 6, 2, 4, 3]$
$g_{37,4} = [8, 7, 6, 5, 4, 3, 2, 4, 5, 6, 7, 8, 1, 3, 4, 5, 6, 2, 4, 5]$
$g_{37,5} = [8, 7, 6, 5, 4, 3, 2, 4, 5, 6, 7, 1, 3, 4, 5, 6, 2, 4, 3, 1]$
$g_{37,6} = [8, 7, 6, 5, 4, 3, 2, 4, 5, 6, 7, 1, 3, 4, 5, 6, 2, 4, 5, 3]$
$g_{38,1} = [8, 7, 6, 5, 4, 3, 2, 4, 5, 6, 7, 8, 1, 3, 4, 5, 6, 7, 2]$
$g_{38,2} = [8, 7, 6, 5, 4, 3, 2, 4, 5, 6, 7, 8, 1, 3, 4, 5, 2, 4, 3]$
$g_{38,3} = [8, 7, 6, 5, 4, 3, 2, 4, 5, 6, 7, 8, 1, 3, 4, 5, 6, 2, 4]$
$g_{38,4} = [8, 7, 6, 5, 4, 3, 2, 4, 5, 6, 7, 1, 3, 4, 5, 2, 4, 3, 1]$
$g_{38,5} = [8, 7, 6, 5, 4, 3, 2, 4, 5, 6, 7, 1, 3, 4, 5, 6, 2, 4, 3]$
$g_{38,6} = [8, 7, 6, 5, 4, 3, 2, 4, 5, 6, 7, 1, 3, 4, 5, 6, 2, 4, 5]$
$g_{39,1} = [8, 7, 6, 5, 4, 3, 2, 4, 5, 6, 7, 8, 1, 3, 4, 5, 6, 7]$
$g_{39,2} = [8, 7, 6, 5, 4, 3, 2, 4, 5, 6, 7, 8, 1, 3, 4, 5, 2, 4]$
$g_{39,3} = [8, 7, 6, 5, 4, 3, 2, 4, 5, 6, 7, 8, 1, 3, 4, 5, 6, 2]$
$g_{39,4} = [8, 7, 6, 5, 4, 3, 2, 4, 5, 6, 1, 3, 4, 5, 2, 4, 3, 1]$
$g_{39,5} = [8, 7, 6, 5, 4, 3, 2, 4, 5, 6, 7, 1, 3, 4, 5, 2, 4, 3]$
$g_{39,6} = [8, 7, 6, 5, 4, 3, 2, 4, 5, 6, 7, 1, 3, 4, 5, 6, 2, 4]$
$g_{40,1} = [8, 7, 6, 5, 4, 3, 2, 4, 5, 6, 7, 8, 1, 3, 4, 5, 2]$
$g_{40,2} = [8, 7, 6, 5, 4, 3, 2, 4, 5, 6, 7, 8, 1, 3, 4, 5, 6]$
$g_{40,3} = [8, 7, 6, 5, 4, 3, 2, 4, 5, 6, 1, 3, 4, 5, 2, 4, 3]$
$g_{40,4} = [8, 7, 6, 5, 4, 3, 2, 4, 5, 6, 7, 1, 3, 4, 5, 2, 4]$
$g_{40,5} = [8, 7, 6, 5, 4, 3, 2, 4, 5, 6, 7, 1, 3, 4, 5, 6, 2]$
$g_{41,1} = [8, 7, 6, 5, 4, 3, 2, 4, 5, 6, 7, 8, 1, 3, 4, 2]$
$g_{41,2} = [8, 7, 6, 5, 4, 3, 2, 4, 5, 6, 7, 8, 1, 3, 4, 5]$
$g_{41,3} = [8, 7, 6, 5, 4, 3, 2, 4, 5, 6, 1, 3, 4, 5, 2, 4]$
$g_{41,4} = [8, 7, 6, 5, 4, 3, 2, 4, 5, 6, 7, 1, 3, 4, 5, 2]$
$g_{41,5} = [8, 7, 6, 5, 4, 3, 2, 4, 5, 6, 7, 1, 3, 4, 5, 6]$
$g_{42,1} = [8, 7, 6, 5, 4, 3, 2, 4, 5, 6, 7, 8, 1, 3, 4]$
$g_{42,2} = [8, 7, 6, 5, 4, 3, 2, 4, 5, 6, 1, 3, 4, 5, 2]$
$g_{42,3} = [8, 7, 6, 5, 4, 3, 2, 4, 5, 6, 7, 1, 3, 4, 2]$
$g_{42,4} = [8, 7, 6, 5, 4, 3, 2, 4, 5, 6, 7, 1, 3, 4, 5]$
$g_{43,1} = [8, 7, 6, 5, 4, 3, 2, 4, 5, 6, 7, 8, 1, 3]$
$g_{43,2} = [8, 7, 6, 5, 4, 3, 2, 4, 5, 6, 1, 3, 4, 2]$
$g_{43,3} = [8, 7, 6, 5, 4, 3, 2, 4, 5, 6, 1, 3, 4, 5]$
$g_{43,4} = [8, 7, 6, 5, 4, 3, 2, 4, 5, 6, 7, 1, 3, 4]$
$g_{44,1} = [8, 7, 6, 5, 4, 3, 2, 4, 5, 6, 7, 8, 1]$
$g_{44,2} = [8, 7, 6, 5, 4, 3, 2, 4, 5, 1, 3, 4, 2]$
$g_{44,3} = [8, 7, 6, 5, 4, 3, 2, 4, 5, 6, 1, 3, 4]$
$g_{44,4} = [8, 7, 6, 5, 4, 3, 2, 4, 5, 6, 7, 1, 3]$
$g_{45,1} = [8, 7, 6, 5, 4, 3, 2, 4, 5, 6, 7, 8]$
$g_{45,2} = [8, 7, 6, 5, 4, 3, 2, 4, 5, 1, 3, 4]$
$g_{45,3} = [8, 7, 6, 5, 4, 3, 2, 4, 5, 6, 1, 3]$
$g_{45,4} = [8, 7, 6, 5, 4, 3, 2, 4, 5, 6, 7, 1]$
$g_{46,1} = [8, 7, 6, 5, 4, 3, 2, 4, 5, 1, 3]$
$g_{46,2} = [8, 7, 6, 5, 4, 3, 2, 4, 5, 6, 1]$
$g_{46,3} = [8, 7, 6, 5, 4, 3, 2, 4, 5, 6, 7]$
$g_{47,1} = [8, 7, 6, 5, 4, 3, 2, 4, 1, 3]$
$g_{47,2} = [8, 7, 6, 5, 4, 3, 2, 4, 5, 1]$
$g_{47,3} = [8, 7, 6, 5, 4, 3, 2, 4, 5, 6]$
$g_{48,1} = [8, 7, 6, 5, 4, 3, 2, 4, 1]$
$g_{48,2} = [8, 7, 6, 5, 4, 3, 2, 4, 5]$
$g_{49,1} = [8, 7, 6, 5, 4, 3, 2, 1]$
$g_{49,2} = [8, 7, 6, 5, 4, 3, 2, 4]$
$g_{50,1} = [8, 7, 6, 5, 4, 3, 1]$
$g_{50,2} = [8, 7, 6, 5, 4, 3, 2]$
$g_{51,1} = [8, 7, 6, 5, 4, 2]$
$g_{51,2} = [8, 7, 6, 5, 4, 3]$
$g_{52,1} = [8, 7, 6, 5, 4]$
$g_{53,1} = [8, 7, 6, 5]$
$g_{54,1} = [8, 7, 6]$
$g_{55,1} = [8, 7]$
$g_{56,1} = [8]$
$g_{57,1} = [\,]$

17